\documentclass[11pt,a4paper]{article}
\usepackage[top=2.5cm,bottom=2.5cm,left=2.5cm,right=2.5cm]{geometry}
\usepackage[utf8]{inputenc}
\usepackage{xcolor}
\usepackage{graphicx}
\usepackage{amsmath,amssymb,latexsym}
\usepackage{epsf,enumerate}
\usepackage{tcolorbox}
\usepackage{url}

\usepackage{bm}

\newtheorem{problem}{Problem}[section]
\newtheorem{conjecture}{Conjecture}[section]

\newenvironment{boxedproblem}{\begin{problem}}{\end{problem}}
\newenvironment{boxedconjecture}{\begin{conjecture}}{\end{conjecture}}

\tcolorboxenvironment{boxedproblem}{}
\tcolorboxenvironment{boxedconjecture}{}

\usepackage[super]{nth}

\usepackage{hyperref}

\makeatletter
\let\Hy@linktoc\Hy@linktoc@none
\makeatother

\title{Open problems of the \nth{32} Workshop on Cycles and Colourings}

\author{(edited by Alfréd Onderko)}

\date{}

\begin{document}

\maketitle

\begin{abstract}
Since its beginnings, every Cycles and Colourings workshop holds one or two open problem sessions; this document contains the problems (together with notes regarding the current state of the art and related bibliography) presented by participants of the \nth{32} edition of the workshop which took place in Poprad, Slovakia during September 8 -- 13, 2024 (see the workshop webpage \url{https://candc.upjs.sk}).
    
\end{abstract}

\tableofcontents

\newpage
\section[Crumby colorings (János Barát)]{Crumby colorings}

\begin{flushleft}
\textsc{\Large{János Barát}}\smallskip \\ Alfréd Rényi Institute of Mathematics, Budapest, Hungary \\ University of Pannonia, Veszprém, Hungary \\ {\tt barat@renyi.hu}
\end{flushleft}

% \begin{center}
% \textsc{\Large{János Barát}}\smallskip \\ Alfréd Rényi Institute of Mathematics, Budapest, Hungary \\ University of Pannonia, Veszprém, Hungary \\ {\tt barat@renyi.hu}
% \end{center}
\bigskip
Consider a blue-red vertex coloring of a 3-connected cubic graph $G$. 
Denote by $B$ and $R$ the subgraphs of $G$ induced on blue and red vertices, respectively.
In~\cite{thomassen2018square} Thomassen conjectured that there is always such a blue-red coloring of $G$ where $\Delta(B) \leq 1$, $\delta(R) \geq 1$ and $R$ does not contain a path with 3 edges.
% In~\cite{} Thomassen conjectured that every 3-connected cubic graph admits a vertex coloring using blue and red such that the maximum degree of a subgraph induced on blue vertices is at most 1, and the minimum degree of the subgraph induced on red vertices is at least 1 and it does not contain a 3-edge path.
Such a coloring is known as a crumby coloring.
If Thomassen's conjecture on crumby colorings were true, it would imply that the square of every planar cubic graph is 7-colorable (Wegner's conjecture proposed in~\cite{wegner1977graphs} and proved in~\cite{thomassen2018square}). However, the conjecture on crumby colorings is false. 
Hence, the following question arises:

\begin{boxedproblem}
    What relaxation of the conditions on $R$ guarantees that such a coloring of a 3-connected cubic graph always exists?
\end{boxedproblem}

\nocite{barat2023crumby}

%\bibliographystyle{plain}
%\bibliography{bibliographies/bib_barat}
%\putbib[bibliographies/bib_barat]
%\putbib

\newpage
\section[Various open problems related to list packing (Stijn Cambie)]{Various open problems related to list packing}

\begin{flushleft}
\textsc{\Large{Stijn Cambie}}\smallskip \\
KU Leuven Kulak, Kortrijk, Belgium\\
{\tt stijn.cambie@kuleuven.be}
\end{flushleft}

\bigskip
A list-assignment $L$ of a graph $G$ is a function $L:V(G)\to 2^{\mathbb{N}}$. If $L:V(G) \rightarrow \binom{\mathbb{N}}{k},$ it is a $k$-list-assignment of $G$.

An $L$-colouring of $G$ is a colouring $c \colon V(G)\to \mathbb{N}$ such that for every vertex $v$, $c(v)\in L(v)$. 
It is proper if adjacent $u,v$ satisfy $c(u)\not=c(v).$

Two $L$-colourings $c_i, c_j$ are disjoint if $c_i(u) \not= c_j(u)$ for every $u \in V(G).$

The {list-chromatic number} (or {choosability}) $\chi_{\ell}(G)$ of a graph is the minimum integer $k$ such that every $k$-list-assignment $L$ admits an $L$-colouring.
Given a $k$-list-assignment $L$ of $G$ we call a collection of $k$ pairwise-disjoint $L$-colourings an $L$-packing of size $k$, or less specifically a list-packing.

The {list packing number} $\chi^\star_\ell(G)$ of $G$ is the least $k$ such that $G$ admits an $L$-packing of size $k$ for any $k$-list-assignment $L$ of $G$.

\subsection*{The list packing conjecture}

When the list packing number $\chi_{\ell}^\star$ was first introduced in~\cite{CCDK21}, inspired by~\cite{AFH96}, the main conjecture appeared if $\chi_{\ell}^\star$ is linearly bounded by $\chi_{\ell}$.
Meanwhile, these are still open:

\begin{boxedproblem}
    Is there a graph $G$ for which $\chi_{\ell}^\star(G)>\chi_\ell(G)+1$?
\end{boxedproblem}
\begin{boxedproblem}
    Is it true that $\chi_{\ell}^\star(G)\le 2\chi_\ell(G)$ for every graph $G$?
\end{boxedproblem}

\subsection*{Planar bipartite graphs}

Alon and Tarsi~\cite{AT92} proved that every planar bipartite graph $G$ satisfies $\chi_{\ell}(G) \le 3.$
That is, if every vertex $v \in V(G)$ has three different possible colours, from a list $L(v)$, one can select a colour $c(v) \in L(v)$ for every $v \in V(G)$ such that no two neighbours $u,v$ satisfy $c(u)=c(v).$
For a $3$-list-assignemnt of a cubic planar graph, we have the following problems:

\begin{boxedproblem}
    Is it possible to choose two disjoint proper colourings?
\end{boxedproblem}
\begin{boxedproblem}
    Is it possible to choose three disjoint proper colourings?
\end{boxedproblem}

In the case of cubic bipartite graphs, it is possible to find two disjoint proper colourings (problem similar to Problem 3) [unpublished].
Hence, the following question arises:
Is it possible to find three disjoint proper colourings in the case of cubic bipartite graphs?

\subsection*{Upper bound by a function of maximum degree}

In~\cite{Cambie_Cames_van_Batenburg_Davies_Kang_2024}, the upper bound for $\chi_{\ell}^\star(G)$ among graphs with maximum degree bounded by $3$ was studied.
It is observed that the analogue of Brook's theorem is also different.
The main question in this paper is the following one:

\begin{boxedproblem}
    Is it true that $\chi_{\ell}^\star(G) \le \Delta(G)+1$ for every graph $G$?
\end{boxedproblem}

\subsection*{Planar graphs}

Almost simultaneously (even while in a different year),~\cite{CCZ23} and~\cite{CS-R24+} gave improved upper bounds for the list packing numbers of planar graphs. Is it possible to improve these bounds?

\begin{boxedproblem}
    Is there a planar graph with $\chi_{\ell}^\star(G) \ge 6$?
\end{boxedproblem}

\begin{boxedproblem}
    Is it true that $\chi_{\ell}^\star(G) \le 7$ for every planar graph?
\end{boxedproblem}

\begin{boxedproblem}
    For a triangle-free planar graph $G$, is it true that $\chi_{\ell}^\star(G) \le 4$? 
\end{boxedproblem}

\subsection*{Removing an edge or vertex}

For most chromatic numbers, it is trivial that removing a vertex or edge cannot change the chromatic number by more than one.
This is not true, for the list packing number.
The following questions posed in~\cite{CH23}, are nevertheless not trivial.

\begin{boxedproblem}
    Let $v$ and $e$ be a vertex and an edge of a graph $G$.
    Is it true that
    $\chi_{\ell}^\star(G) \leq \chi_{\ell}^\star(G \setminus v)+2$?

    \medskip
    Or, is it true that $\chi_{\ell}^\star(G) \leq \chi_{\ell}^\star(G \setminus e)+2$?
\end{boxedproblem}

\subsection*{Computational hardness}

Last, we note that little is known about the computation hardness of list packing.
In~\cite{2024arXiv240203520C}, there is a conjecture on a related question.

\begin{boxedproblem}
    Fix $k>2$. Given a graph $G$ as input, is the problem of deciding whether $\chi^\star_\ell(G) \le k$ complete for the complexity class $\Pi^{\bf P}_2$?
\end{boxedproblem}

\begin{boxedproblem}[Question 7.3 in~\cite{CC24}]
    Fix $t>1$. 
    Is the problem of determining the value of $\chi^\star_\ell(G)$ or $\chi^\bullet_\ell(G)$, given a graph $G$ of treewidth at most $t$ as input, in $P$?
\end{boxedproblem}

\newpage
\section[Chromatic index of extended Petersen graphs (Ge\v na Hahn)]{Chromatic index of extended Petersen graphs}

\begin{flushleft}
    \textsc{\Large{Ge\v na Hahn }}\smallskip\\
    Universit\'e de Montr\'eal, Montr\'eal, Canada\\
    {\tt hahn@iro.umontreal.ca}
\end{flushleft}

\bigskip
Let $n = 4k + 1$ for some integer $k \geq 1$, and let $S \subseteq \{1,2, \dots, 2k\}$ such that $|S| = k$. The extended Petersen graph $\mathrm{XP}(2n,S)$ is a graph with the vertex set $V = \mathbb{Z}_n \times \{1,2\}$, and the edge set 
\begin{align*}
    E = \{x_1y_1 \colon d(x_1,y_1) \in S\} \cup \{x_2y_2 \colon d(x_2,y_2) \in \overline{S} \} \cup \{x_1x_2 \colon x \in \mathbb{Z}_n\},
\end{align*}
where $\overline{S} = \{1, \dots, 2k \} \setminus S$ and $d(x,y) = \min\{|x-y| , n - |x-y|\}$.

Extended Petersen graphs were introduced in~\cite{HorakRosa} by Hor\'ak and Rosa, who proposed the following conjecture:
\begin{boxedconjecture}[Horák, Rosa~\cite{HorakRosa}]
    Every extended Petersen graph, except Petersen graph itself, is of Class 1.
\end{boxedconjecture}

% \bibliographystyle{plain}
% \bibliography{bibliographies/bib_hahn_1}
%\putbib[bibliographies/bib_hahn_1]
%\putbib

\newpage
\section[A conjecture on rainbow connectivity (Ge\v na Hahn)]{A conjecture on rainbow connectivity}

\begin{flushleft}
    \textsc{\Large{Ge\v na Hahn }}\smallskip\\
    Universit\'e de Montr\'eal, Montr\'eal, Canada\\
    {\tt hahn@iro.umontreal.ca}
\end{flushleft}

\bigskip
In~\cite{chartrand2016highly} Chartrand proposed a concept of rainbow connectivity: an edge-colored graph $G$ is rainbow connected if every two vertices of $G$ are connected by a rainbow path (that is a path with all its edges colored with different colors). 
The minimum number of colors required so that the colored graph is rainbow connected, is called the rainbow connection number, and it is denoted by $\mathrm{rc}(G)$.

\begin{boxedproblem}
     Does there exist a constant $c \geq 0$ such that if $G$ is a non-complete graph of order $n$ with $\delta(G) \geq \frac{n}{2}+c$, then $\mathrm{rc}(G) = 2$? 

    \medskip
     In particular, does $\delta(G) \geq \frac{n}{2}$ imply that $\mathrm{rc}(G) = 2$?
\end{boxedproblem}

It is known due to Gimbel~\cite{gimbel} that if $\delta(G) \geq \frac{n}{2} + \log_2(n) - 1$ then $\mathrm{rc}(G) = 2$.

% \bibliographystyle{plain}
% \bibliography{bibliographies/bib_hahn2}
%\putbib[bibliographies/bib_hahn2]

\newpage
\section[Edge-colorability and perfect matchings of $5$-edge-connected $5$-regular graphs (Davide Mattiolo)]{Edge-colorability and perfect matchings of $5$-edge-connected $5$-regular graphs}

\begin{flushleft}
    \textsc{\Large{Davide Mattiolo}}\footnote{Supported by a Postdoctoral Fellowship of the Research Foundation Flanders (FWO), grant 1268323N.}\smallskip\\
    (joint work with Y.\ Ma, E.\ Steffen and I.\ H.\ Wolf)\\
    KU Leuven Kulak, Kortrijk, Belgium\\
    {\tt davide.mattiolo@kuleuven.be}
\end{flushleft}

\bigskip

    A graph $G$ is {Class} $1$ if it is $\Delta(G)$-edge-colorable, otherwise it is {Class} $2$.
    
    It is known that $r$-edge-connected $r$-regular Class $2$ graphs exist for all integers $r\ge2$, with $r\ne5$ (see \cite{Meredith} and \cite{MMSW}).
    
    Surprisingly, it seems that no $5$-edge-connected $5$-regular Class $2$ graph is known.
    
    \begin{boxedproblem}[Ma, Mattiolo, Steffen, Wolf~\cite{MMSW}]
    Is there any $5$-edge-connected $5$-regular Class $2$ graph?
    \end{boxedproblem}
    
    A negative answer to Problem~1 would imply the Berge-Fulkerson and the $5$-cycle double cover conjectures, see \cite{MMSW} for more details.
    
    An $r$-{graph} is an $r$-regular graph $G$ such that every odd set $X \subseteq V(G)$ is connected by at least $r$ edges to its complement $V(G)\setminus X$.
    An $r$-graph is {poorly matchable} if any two of its perfect matchings intersect. In \cite{Rizzi} it is shown that there are poorly matchable $r$-graphs for every $r\ge3.$ All such constructed graphs contain a $4$-edge-cut. Moreover, Thomassen \cite{THOMASSEN2020343} conjectured that there is a natural number $r_0$ such that, for every $r\ge r_0$, every $r$-edge-connected $r$-graph has two disjoint perfect matchings. Since the poorly matchable $4$-graphs constructed in \cite{Rizzi} are $4$-edge-connected, such an integer $r_0$ must be at least $5$.
    
    \begin{boxedproblem}
        Is there a poorly matchable $5$-edge-connected $5$-graph?
    \end{boxedproblem}
    
    In \cite{MMSW} it is proved that, if every $5$-edge-connected $5$-graph has two edge-disjoint perfect matchings then the Fan-Raspaud Conjecture holds.

% \bibliographystyle{plain}
% \bibliography{bibliographies/bib_mattiolo}

\begin{thebibliography}{1}

\bibitem{barat2023crumby}
Barát, J., Blázsik, Z. L., Damásdi, G.
\newblock Crumby colorings---Red-blue vertex partition of subcubic graphs regarding a conjecture of Thomassen.
\newblock {\em Discrete Mathematics} 346(4) (2023), 113281.

\bibitem{thomassen2018square}
Thomassen, C.
\newblock The square of a planar cubic graph is 7-colorable.
\newblock {\em Journal of Combinatorial Theory, Series B}, 128 (2018), 192--218.

\bibitem{wegner1977graphs}
Wegner, G.
\newblock Graphs with given diameter and a coloring problem.
\newblock {\em Technical report}, (1977).

\end{thebibliography}

\begin{thebibliography}{1}

\bibitem{AT92}
Alon, N., Tarsi, M.
\newblock Colorings and orientations of graphs.
\newblock {\em Combinatorica}, 12(2) (1992), 125--134.

\bibitem{AFH96}
Alon, N., Fellows, M. R., Hare, D. R.
\newblock Vertex transversals that dominate.
\newblock {\em Journal of Graph Theory}, 21(1) (1996), 21--31.

\bibitem{CC24}
Cambie, S., van Batenburg, W. C.
\newblock {Fractional list packing for layered graphs}.
\newblock {\em arXiv e-prints}, (2024), arXiv:2410.02695.

\bibitem{CCDK21}
Cambie, S., van Batenburg, W. C., Davies, E., Kang, R. J.
\newblock Packing list-colorings.
\newblock {\em Random Structures \& Algorithms}, 64(1) (2024), 62--93.

\bibitem{Cambie_Cames_van_Batenburg_Davies_Kang_2024}
Cambie, S., van Batenburg, W. C., Davies, J., Kang, R. J.
\newblock List packing number of bounded degree graphs.
\newblock {\em Combinatorics, Probability and Computing}, 33(6) (2024), 807--828.

\bibitem{CCZ23}
Cambie, S., van Batenburg, W. C., Zhu, X.
\newblock {Disjoint list-colorings for planar graphs}.
\newblock {\em arXiv e-prints}, (2023), arXiv:2312.17233.

\bibitem{CH23}
Cambie, S., {H{\"a}m{\"a}l{\"a}inen}, R.
\newblock {Packing colourings in complete bipartite graphs and the inverse problem for correspondence packing}.
\newblock {\em arXiv e-prints}, (2023), arXiv:2303.01944.

\bibitem{2024arXiv240203520C}
Camrud, E., Davies, E., Karduna, A., Lee, H.
\newblock {Sampling List Packings}.
\newblock {\em arXiv e-prints}, (2024), arXiv:2402.03520.

\bibitem{CS-R24+}
Cranston, D. W., Smith-Roberge, E.
\newblock {List-Coloring Packing and Correspondence-Coloring Packing of Planar Graphs}.
\newblock {\em arXiv e-prints}, (2024), arXiv:2401.01332.

\end{thebibliography}

\begin{thebibliography}{1}

\bibitem{HorakRosa}
Horak, P., Rosa, A.
\newblock Extended Petersen graphs.
\newblock {\em Discrete mathematics}, 299(1-3) (2005), 129--140.

\end{thebibliography}

\begin{thebibliography}{1}

\bibitem{chartrand2016highly}
Chartrand, G.
\newblock Highly irregular.
\newblock {\em Graph Theory: Favorite Conjectures and Open Problems-1}, (2016), 1--16.

\bibitem{gimbel}
Gimbel, J.
\newblock A conjecture on rainbow connectivity.
\newblock https://iti.mff.cuni.cz/series/2024/ 687.pdf.

\end{thebibliography}

\begin{thebibliography}{1}

\bibitem{MMSW}
Ma, Y., Mattiolo, D., Steffen, E., Wolf, I. H.
\newblock Pairwise disjoint perfect matchings in $r$-edge-connected $r$-regular graphs.
\newblock {\em SIAM Journal on Discrete Mathematics}, 37(3) (2023),1548--1565.

\bibitem{Meredith}
Meredith, G. H. J.
\newblock Regular $n$-valent $n$-connected nonhamiltonian non-$n$-edge-colorable graphs.
\newblock {\em Journal of Combinatorial Theory, Series B}, 14(1) (1973),55--60.

\bibitem{Rizzi}
Rizzi, R.
\newblock Indecomposable $r$-graphs and some other counterexamples.
\newblock {\em Journal of Graph Theory}, 32(1) (1999), 1--15.

\bibitem{THOMASSEN2020343}
Thomassen, C.
\newblock Factorizing regular graphs.
\newblock {\em Journal of Combinatorial Theory, Series B}, 141 (2020), 343--351.

\end{thebibliography}

\begin{thebibliography}{1}

\bibitem{CT-ar}
Caro, Y., Tuza, Z.
\newblock {Monochromatic graph decompositions inspired by anti-Ramsey theory and the odd-coloring problem}.
\newblock {\em arXiv e-prints}, (2024), arXiv:2408.04257.

\end{thebibliography}
%\putbib[bibliographies/bib_mattiolo]

\newpage
\section[Subcubic $K_3$-free graphs of Class 2 (Ingo Schiermeyer)]{Subcubic $K_3$-free graphs of Class 2}

\begin{flushleft}
    \textsc{\Large{Ingo Schiermeyer}}\smallskip\\
    Technische Universität Bergakademie Freiberg, Freiberg, Germany\\
    AGH University of Krakow, Krakow, Poland\\
    {\tt Ingo.Schiermeyer@math.tu-freiberg.de}
\end{flushleft}

\bigskip
It is known that subdividing a single edge of a $3$-edge colorable cubic graph yields an overfull subcubic graph $G$. Hence, $G$ is of Class 2, i.e., $\chi'(G) = 4$. 
Starting with a $K_3$-free cubic graph, we produce, in that way, a $K_3$-free subcubic graph of Class 2. 
This raises the following question:
\begin{boxedproblem}
    Is it possible to describe all subcubic $K_3$-free graphs of Class 2?
\end{boxedproblem}
This class of graphs contains all snarks, which are cubic 
graphs of Class 2. The existence of snarks is a challenging topic on its own.

\newpage
\section[Special orderings of vertices (Zsolt Tuza)]{Special orderings of vertices}

\begin{flushleft}
    \textsc{\Large{Zsolt Tuza}}\smallskip\\
    University of Pannonia, Veszprém, Hungary\\
    Alfréd Rényi Institute of Mathematics, Budapest, Hungary\\
    {\tt tuza.zsolt@mik.uni-pannon.hu}
\end{flushleft}

\bigskip
The following questions are raised in a joint paper with Yair Caro \cite{CT-ar}, concerning some variants of anti-Ramsey problems.

Suppose that the vertices of a graph $G$ are put in a linear order $v_1,\dots,v_n$.
For $i=1,\dots,n$ define $d^-(v_i) = |\{j < i \colon v_iv_j \in E(G)\}|$ and $d^+(v_i) = |\{j > i \colon v_iv_j \in E(G)\}|$.

We are interested in the existence of two types of orderings (and a way to find them as efficiently as possible, if they exist) given by the following two conditions:
\begin{itemize}
    \item Type-A order: for each $v \in V(G)$, $d^-(v)$ is odd or equal to 0; 
    \item Type-B order: for each $v \in V(G)$, $d^-(v)$ is odd or equal to 0 (as in a Type-A order), or $d^-(v)$ is even and $d^+(v) = 0$.
\end{itemize}
These types of orders are useful in the context of anti-Ramsey problems and odd-coloring \cite{CT-ar}.
We propose the following algorithmic questions.

\begin{boxedproblem}
    Can it be decided in polynomial time whether a graph $G$ admits a Type-A order, or a Type-B order?
\end{boxedproblem}

\begin{boxedproblem}
    If a Type-A or Type-B order exists, can it be found in polynomial time?
\end{boxedproblem}

Of course, for a given linear order of the vertices, checking whether it is a Type-A or Type-B order can be tested in linear time.

\end{document}